\documentclass[11pt,a4paper]{amsart}

\usepackage[T1]{fontenc}
\usepackage[utf8]{inputenc}
\usepackage[a4paper,margin=1in]{geometry}
\usepackage{lmodern}
\usepackage{microtype}
\usepackage{amssymb,amsmath,amsthm}
\usepackage{graphicx,enumerate}
\usepackage[mathlines]{lineno}
\usepackage[hidelinks,bookmarks=false]{hyperref}

\newcommand\cG{{\mathcal G}}
\newcommand\cH{{\mathcal H}}

\newcommand{\ex}{\mathop{}\!\mathrm{ex}}

\makeatletter
\newtheorem*{rep@theorem}{\rep@title}
\newcommand{\newreptheorem}[2]{%
\newenvironment{rep#1}[1]{%
 \def\rep@title{#2 \ref{##1}}%
 \begin{rep@theorem}}%
 {\end{rep@theorem}}}
\makeatother

\newtheorem{thm}{Theorem}%[chapter]
\newtheorem{theorem}{Theorem}

\newtheorem{lemma}[thm]{Lemma}

\newtheorem{proposition}[thm]{Proposition}

\newtheorem{problem}[thm]{Problem}

\newcommand\cref[1]{Corollary~\ref{cor:#1}}

\title{Edge mappings of graphs: Ramsey-type parameters}

\author{Yair Caro}
\address{Department of Mathematics, University of Haifa-Oranim, Israel}
\email{yacaro@kvgeva.org.il}

\author{Bal\'azs Patk\'os}
\address{HUN-REN Alfr\'ed R\'enyi Institute of Mathematics}
\email{patkos@renyi.hu}
\thanks{Research of Patk\'os, Tuza, and Vizer is partially supported by NKFIH grants SNN 129364 and FK 132060.}

\author{Zsolt Tuza}
\address{HUN-REN Alfr\'ed R\'enyi Institute of Mathematics and University of Pannonia}
\email{tuza.zsolt@mik.uni-pannon.hu}

\author{M\'at\'e Vizer}
\address{Department of Computer Science and Information Theory, Budapest University of Technology and Economics}
\email{vizermate@gmail.com}


\date{}

\begin{document}
%\linenumbers

\begin{abstract}
    In this paper, we launch a systematic study of parameters concerning edge mappings of graphs. Inspired by Ramsey's theorem, the quantity $m(G,H)$ is defined to be the smallest integer $n$ such that for every $f:E(K_n)\rightarrow E(K_n)$ either there is an $f$-fixed copy of $G$ with $f(e)=e$ for all $e\in E(G)$, or an $f$-free copy of $H$ with $f(e)\notin E(H)$ for all $e\in E(H)$. Incorporating new ideas, we extend many old results from the 1980s and prove many new exact results, mostly concerning $m(T,K_r)$, where $T$ is a tree, and $m(G,K_{1,r})$. We also study further related parameters, most of them introduced in the 1980s, and obtain substantial progress regarding these parameters.
\end{abstract}

\maketitle

\section{Introduction}

In a recent paper \cite{CPTV}, we investigated Turán-type problems related to edge mappings in various settings. The inception of our investigation found its roots in the systematic inquiries conducted by Erd\H os and Hajnal \cite{EH58}, who initiated the study of edge-mappings (alternatively called set-mappings). Their primary focus was directed towards infinite sets; however, their exploration also led to the initial examination of such problems in the finite case, influenced by Ramsey's theorem. In this paper, our study delves deeper into the interplay between edge mappings and Ramsey-type problems.

\bigskip

There are several very interesting classes of edge mappings and graph invariants
related to them. In order to present our results, in the next subsection we
introduce those definitions and necessary notation.

\subsection{Definitions and notation}
Let $E_n  = E(K_n)$  denote the set of edges of the complete graph $K_n$ on $n$ vertices, and let $f: E_n \rightarrow  E_n$ be a function. A subgraph $G$ of $K_n$ is called 

\begin{itemize}

\smallskip

\item \textit{$f$-fixed,} if $f( e) = e$ for all $e  \in E(G)$,
\smallskip

\item \textit{$f$-shifted,} if $f( e)\neq  e$ for all $e \in E(G)$,
\smallskip

\item \textit{$f$-strong shifted,} if $|f(e) \cap  e| = 0$ for all $e \in E(G)$,
\smallskip

\item \textit{$f$-free,} if $f(e) \notin E(G)$ for all $e \in E(G)$, and
\smallskip

\item \textit{$f$-exclusive,} if $|f(e ) \cap   V(G) | = 0$ for all $e\in E(G)$.

\smallskip

\end{itemize}

\medskip

\noindent
Let us introduce some notation for edge mappings. Let $F_n  :=   \{ f : E_n\rightarrow E_n \}$ and $d\in \{0,1\}$.  

\begin{itemize}

\item Let $F_{n,d} := \{  f \in F_n : \ |f(e) \cap  e| \le d  ~\text{for every}\ e \in E_n \}$.

%\smallskip 

%\item For $0\le m \le \binom{n}{2}$ let \[F_{n,m,d}:=\{  f \in F_n  : \ |f(e )\cap e|\le d ~\text{for at least $m$ edges}  \ e \in E_n \}.\]

\smallskip 

\item For a graph $H$ let $$F_{H,d}:=\{f:E(H)\rightarrow E(H) ~\text{with} \  |f(e )\cap e|\le d   ~\text{for every} \ e \in E(H) \} .$$

\end{itemize}

\vskip 0.3truecm

\noindent

\subsection{Parameters related to Ramsey numbers.} Now we introduce parameters related to Ramsey numbers of graphs, which will be the main focus of our investigations. Some of them were defined in \cite{AC,C}, some of them are new. 
\smallskip 

Let us recall that for two graphs $G$ and $H$ the \textit{Ramsey number} $R(G, H)$ is the smallest integer $n$ such that in every red-blue coloring of $E(K_n)$ there exists either a red copy of $G$ or a blue copy of $H$. 

%Also, the \textit{Tur\'an number} $\ex(n,G)$ is the maximum number of edges that an $n$-vertex graph can have without containing $G$ as a subgraph.

\vskip 0.2truecm

For two finite graphs $G$ and $H$, we define
\begin{align*}
z(G,H) := \min\{ n :{}& \text{for every } f\in F_n, \text{ there is either an } f\text{-fixed copy of }G\\
&\text{or an } f\text{-shifted copy of }H\},\\
m(G,H) := \min\{ n :{}& \text{for every } f\in F_n, \text{ there is either an } f\text{-fixed copy of }G\\
&\text{or an } f\text{-free copy of }H\}.
\end{align*}

\medskip 

One may want to define $m^*(G, H)$ as the minimum $n$  such that every  $f \in F_n$  either admits an $f$-fixed copy of $G$ or an $f$-exclusive copy of $H$. The following mapping shows that we need to alter the definition:
let $V(K_n)  = \{1,2,\dots,n\}$ and, for $1\le x<y\le n$, $y\ge 3$, define $f(x,y)=(z,y)$, where $z\equiv x+1$ mod $y$; finally, let $f(1,2)=(2,3)$. Then no edge is $f$-fixed. Furthermore, since $|f(e)\cap e|=1$ holds for all edges $e$, there do not exist any  $f$-exclusive copies of any nontrivial graph $H$. 

\smallskip 

Therefore we exclude the option $|f(e)\cap e|=1$ and introduce the definition in the following form: 
\begin{itemize}
    \item 
    $m^*(G,H) := \min\{n : \forall f \in F_n ~\text{with}\ |e\cap f(e)|\neq 1$ for all $e\in E_n$,  there is either an $f$-fixed copy of $G$ or an $f$-exclusive copy of $H\}$.
\end{itemize}

Further, for a graph $G$, we define

\begin{itemize}

\item $w(G) := \min \{ n : \forall f  \in F_{n,0}, ~\text{there is an $f$-exclusive copy of}\ G ~\text{in}\ K_n \}$.

\end{itemize}

For  $d = 0, 1$ and a graph $G$, we define  

\begin{itemize}

\item $g(G,d) := \min\{n : \forall f\in F_{n,d}, ~\text{there is an $f$-free copy of}\ G  \}$.  
We write $g(G)$ for $g(G,1)$. 

\end{itemize}

\medskip

For any pair $G,H$ of graphs, $G\cup H$ denotes their vertex-disjoint union, $G+ H$  their join, and for a positive integer $k$ we denote by $kG$ the vertex-disjoint union of $k$ copies of $G$.

\bigskip

\subsection{Our results.} 

The Tur\'an number $\ex(n,F)$ is the maximum number of edges that an $n$-vertex $F$-free graph can have. We say that a $k$-vertex tree
$T$ is an \textit{Erd\H os--S\'os tree} (EST for short) if $\ex(n,T)\le \frac{k-2}{2}n$. An old conjecture of Erd\H os and S\'os  \cite{E64} states that all trees satisfy this inequality. A $k$-vertex tree $T$ is a \textit{strong Erd\H os--S\'os tree} if $\ex(n,T) < \frac{k-2}{2}n$ for all $n$ not divisible by $k-1$.

\begin{thm}\label{mtkr}

\ 

\begin{enumerate} 

\item Let $T$  be a tree with $k \ge 2$ vertices. Then $$m(T,K_r)\ge (r-1)(k-1)+1.$$

\item Let $T$  be an EST with $k \ge 2$ vertices. 

\begin{enumerate}[a.]
    \item For any $r \ge 2$ there exists a constant $c_r$ such that $m(T,K_r)\le (r-1)(k-1) +c_r$ holds.
    
    \item 
    $m(T ,K_3) \le  2k + 2$, and  if $T\notin \{ K_2, K_{1,3} ,  P_4 \}$,  then  $m(T,K_3) \le 2k +1$.
    %\item
     %$m(T,K_3)  \ge 2k -1$.
     \item
     If $T \in  \{K_2,K_{1,3}, P_4 \}$,  then  $m(T, K_3)  =  2k +2$.
     \item 
     $3k-2\le m(T,K_4)\le 3k+33$.
\end{enumerate}  

\end{enumerate}

\end{thm}

Extending a result from \cite{AC} that states $m(K_{1,t}, K_{1,r})  = t + 2r -  2$, we prove:

\begin{thm}
\label{th:m-tree}
    If $T$ is a tree on $k$ vertices, then the following inequalities hold.
\begin{enumerate}
    \item 
     If $T$ is an EST, then $m(T,K_{1,r} )\le  k +2r  - 2$ and this bound is sharp if $(k-1) \mid 2(r-1)$.
     \item
     If $k-1$ does not divide $2(r-1)$ and $T$ is a strong EST, then $m(T,K_{1,r}) \le  k+2r  - 3$.
     \item
      $m(T,K_{1,r} )\ge k +2r  - 4$ if  $r > k^2$.
     \item
      $m(T,K_{1,r} ) \ge k +2r  - 3$,  if  $k \equiv 1$ $(\mathrm{mod}~2)$  and $r > k^2$.
\end{enumerate}
\end{thm}

In fact, for many graphs $G$  we exactly determine $m(G,K_{1,r})$  provided that the chromatic number $\chi(G)$ is sufficiently large with respect to $r$; see Section 4.2, Theorems 22 and 23. For the function $m^*$ we obtain the following analog of Theorem \ref{th:m-tree}.

\begin{thm}\label{*tk}\
\begin{enumerate}
    \item
    For every Erd\H os--S\'os tree $T$ on $k\ge 3$  vertices, $m^*(T,K_{1,r})  \le  k+5r  - 5$. 
    \item
    Suppose $T$ is a tree on $k=2\ell+1$ vertices. If $k-1$ divides $2r-2$, then $m^*(T ,K_{1,r}) \ge  k+2r- 2$.
    \end{enumerate}
\end{thm}

For general upper and lower bounds concerning $m^*(G,H)$ we refer to Proposition 9.

\ 

\subsection{Estimates on $g(G,1)$ and $g(G,0)$.}  

Illustrative examples of the type of results we prove are the following upper bounds.

\begin{thm}\label{degreethm}
Let $G$ be a graph with $k$ vertices, $m\ge 2$ edges, and
 degree sequence $d_1,\dots,d_k$.
If $4\binom{m}{2} + (n-7) \sum_{i=1}^k \binom{d_i}{2} < (n-2)(n-3)$,
 then $g(G) \le n$.
\end{thm}

\begin{thm}\label{gGkm}
Let $G$ be a graph with $k$ vertices and $m\ge 2$ edges.
If $4m(k-3) \le n(n-1)(n-2)$, then $g(G,0) \le n$.
\end{thm}

\ 

\subsection{Estimates on $w(G)$} Among several results, we prove the following.

\begin{thm}\label{w}
\begin{enumerate}[(i)]
    \item Every graph $G$ with $k\ge 4$ vertices and $m\ge 1$ edges has
    \[
        w(G) \le 2km-4m+2.
    \]
    \item If $G$ is a graph with $k>2$ vertices, $m>1$ edges and $k,m,n$ satisfy $(1-\frac{2}{k})2m>(n-2)\ln n$, then $w(G)>n$, in particular
    \[
        w(G) > \left(1 - \frac{2}{k} - o(1)\right) \frac{2m}{\ln m} \,.
    \]
\end{enumerate}
\end{thm}

We end this subsection by defining two parameters (used already in \cite{CPTV}) that we will use in some of our proofs and that we will relate to parameters defined above.

\begin{itemize}
    \item For a given graph $G$, $d  = 0 ,1$ and $n \ge 3$, let
\begin{eqnarray}
h(n,d,G) := \max\{ m :& & \hspace*{-1.4em}
\text{$\exists$ $H$, $|V(H)|=n$, $|E(H)|=m$, and $\exists f \in F_{H,d}$} \nonumber \\
& &  \hspace*{-1.4em}
\text{ with no $f$-free copy of $G$ in $H$} \}. \nonumber
\end{eqnarray}

We abbreviate $h(n,1,G)$ to  $h(n,G)$.

\item 

For a given graph $G$, and $n \ge 4$,   let 
\begin{eqnarray}
s(n, G) :=  \max\{ m :& &  \hspace*{-1.4em}
\text{$\exists$ $H$, $|V(H)|=n$, $|E(H)|=m$, and $\exists f\in F_{H,0}$ } \nonumber \\
& &  \hspace*{-1.4em}
\text{with no $f$-exclusive  copy of $G$ in $H$} \}. \nonumber
\end{eqnarray}
\end{itemize}

\subsection{Organization of the paper}

Section~2 outlines basic connections among some of the previously defined parameters. Section~3 collects the tools used throughout the paper. In Section~4, we present our contributions concerning $m(G,K_r)$, $m(T,2K_2)$, and $m(G,K_{1,r})$, proving Theorems~\ref{mtkr}, \ref{th:m-tree}, and \ref{*tk}, and further exact results concentrated in Theorems~\ref{properties}, \ref{th:exact-families}, and \ref{th:exact-subgraphs}. In Section~5, we present our contributions regarding $g(G)$ and $g(G,0)$, proving Theorems~\ref{degreethm} and \ref{gGkm}, as well as the exact determination of $g(tK_2,1)$ and $g(tK_2,0)$. Section~6 is devoted to the proof of Theorem~\ref{w}, a further sharper result concerning $w(K_k)$, and an interesting result concerning exact values of $w(K_{1,r})$. The paper concludes in Section~7, where we present a diverse array of open problems for future research.

\section{Basic connections between parameters}

In this section, we gather some simple connections among the parameters defined in the introduction and their relation to the Ramsey and Tur\'an numbers. We first relate the functions $m$ and $z$ to the Ramsey function.

\begin{proposition}\label{mzR}
For any two graphs $G,H$ we have
$m(G,H)\ge z(G,H) = R(G,H)$.  
\end{proposition}

\begin{proof}

The inequality $m(G, H) \ge z(G, H)$ holds by definition, so we only need to prove $z(G, H) = R(G, H)$. 

 Suppose $n\ge R(G,H)$ and $f \in F_n$.   Color $e$ red  if $f( e)  = e$ and color $e$ blue  if $f( e) \neq e$. 
There is either a red hence $f$-fixed copy of $G$ or a blue hence $f$-shifted copy of $H$. This proves $z(G,H)\le R(G,H)$.

Suppose next $n = R(G,H) -1$ and let $c : E(K_n) \rightarrow  \{ red, blue \}$ be a coloring with neither red $G$ nor blue $H$. Define $f : E(K_n)   \rightarrow E(K_n)$ by $f( e)  = e$  if $c( e)$ is red  and $f(e )\neq e$ (just by choosing any other edge for $f( e)$) otherwise. If there is a copy of $G$ which is $f$-fixed then it is monochromatic red which is not the case, and if there is a copy of $f$-shifted $H$ then this copy is monochromatic blue which is not the case either. 
\end{proof}

We continue with an inequality relating the functions $h$, $s$ and $m$ to the Turán function.

\begin{proposition}
    For any two families of graphs $\cG, \cH$ the following bounds are valid:
    \begin{enumerate}
        \item 
        if $\ex(n,\cG) + h(n,\cH) < \binom{n}{2}$, then $m(\cG, \cH)\le n$,
        \item
       if $\ex(n,\cG) + s(n,\cH)  <    \binom{n}{2}$,  then $m^*(\cG,\cH) \le n$. 
    \end{enumerate}
\end{proposition}

\begin{proof}
     If for $f\in F_n$, the number of $f$-shifted edges is greater than $h(n,\cH)$, then by definition of $h$, $f$ admits an $f$-free $H$ for some $H \in \cH$. Otherwise the number of $f$-fixed edges is at least $\binom{n}{2}-h(n,\cH)>\ex(n,\cG)$, so $f$ admits an $f$-fixed $G$ for some $G\in \cG$. This proves (1); the proof of (2) is analogous.
\end{proof}

\begin{proposition}
For any two graphs $G,H$ we have
$R(G,H)\le m^*(G,H)  \le R(G,K_{w(H)})$.
\end{proposition}

\begin{proof}
To see the first inequality, along the lines of Proposition \ref{mzR}, suppose  $n = R(G, H) -1$ and let $c: E(K_n) \rightarrow  \{ red, blue \}$ be a coloring with neither red $G$ nor blue $H$. Define $f : E(K_n)   \rightarrow E(K_n)$ by $f( e)  = e$  if $c( e)$ is red  and $f(e )$ to be any edge disjoint with $e$ if $c(e)$ is blue. Then any $f$-fixed copy of $G$ would be monochromatic red, and any $f$-exclusive copy of $H$ would be monochromatic blue. But none of them exist.
    
For the second inequality, let $n \ge R(G, K_{w(H)})$, and consider a mapping $f\in F_n$ where every edge is either fixed or strong-shifted (observe that $n \ge R(G, K_{w(H)})$ implies $n\ge 4$ which is needed for the existence of such a mapping).
We color an edge red if it is fixed and blue if it is strong-shifted and obtain either a red hence $f$-fixed copy of $G$ or a blue copy of $K_{w(H)}$ that is a copy in which all edges are strong-shifted.  In the former case we are done. In the latter case, we can define $f':E(K_{w(H)})\rightarrow E(K_{w(H)})$ by distinguishing cases according to $|V(K_{w(H)})\cap f(e)|$. We let $f'(e)=f(e)$ if $f(e)\subseteq V(K_{w(H)})$ and  $f'(e)=e^*$ for any $e^*\in E(K_{w(H)})$ with $e\cap e^*=\emptyset$ if $f(e)\cap V(K_{w(H)})=\emptyset$. Finally if $f(e)=xy$ with $f(e)\cap V(K_{w(H)})=\{x\}$, then we let $f'(e)=xz$ for some $z\in V(K_{w(H)})\setminus (e\cup \{x\})$. The function $f'$ is in
$F_{w(H),0}$ and, by the definition of $w(H)$, there must exist an $f'$-exclusive copy of $H$, and clearly, this copy is also $f$-exclusive.
\end{proof}

\section{Tools}

This section gathers tools used several times in the sequel. We start with a simple lemma. Part (1) is from \cite{AC} (Lemma 2.3) and Part (2) already appears in our parallel manuscript \cite{CPTV}, but as we will use them multiple times in later proofs, we restate and reprove them in order to be self-contained.
The \textit{undirected version} of a directed graph $\overrightarrow{\Gamma}=(V,\overrightarrow{E})$ is the undirected graph $G=(V,E)$ with $E=\{uv: \overrightarrow{uv}\in \overrightarrow{E} ~\text{or}\ \overrightarrow{vu}\in \overrightarrow{E}\}$.

\begin{lemma}\label{outdeg}
    Let $\overrightarrow{\Gamma}$ be a directed graph on $n$ vertices with maximum out-degree at most $d$ and let $\Gamma$ be the undirected version of $\overrightarrow{\Gamma}$.
\begin{enumerate}
    \item
    We have $\chi(\Gamma)\le 2d+1$, and thus $\alpha(\Gamma)\ge \frac{n}{2d+1}$.
    \item
    If $d=1$, and $m$ vertices have out-degree $0$, then  we have $\alpha(\Gamma)\ge m+\lceil \frac{n-2m}{3}\rceil$. Furthermore, if $m=0$ and $3|n$, then equality holds if and only if $\overrightarrow{\Gamma}$ consists of vertex-disjoint cyclically directed triangles.
    \end{enumerate}
\end{lemma}

\begin{proof}
    We start with the proof of (2). Note that as the out-degree of each vertex of  $\overrightarrow{\Gamma}$ is at most 1, the components of $\Gamma$ are trees or unicyclic graphs \footnote{A \textit{unicyclic} graph is a connected graph containing exactly one cycle.}. Also, note that the components containing a vertex with 0 out-degree (in $\overrightarrow{\Gamma}$) are trees. Trees are bipartite, unicyclic graphs are 3-colorable. This shows the inequality $\alpha(\Gamma) \ge m + \lceil
\frac{n-2m}{3}\rceil$. 
    
    If $m=0$, then we need all components of $\Gamma$ to be unicyclic. As all such graphs can be three-colored with one color class having size one, we need $\left\lceil \frac{|C|-1}{2} \right\rceil=\frac{|C|}{3}$ for all components $C$. This implies $|C|=3$ and the out-degree condition yields that the triangles are cyclically directed.

\smallskip    

    To see (1), observe that if the out-degree of $\overrightarrow{\Gamma}$ is $d$, then the average degree, and thus the degeneracy\footnote{ A \textit{$k$-degenerate} graph is an undirected graph in which every subgraph has a vertex of degree at most $k$. The \textit{degeneracy} of a graph is the smallest value of $k$ for which it is $k$-degenerate.} of $\Gamma$ is at most $2d$, so its chromatic number is at most $2d+1$.
\end{proof}

%We will apply Lemma \ref{outdeg} in the following situation: given a graph $G$ a mapping $f\in F_{G,1}$, and a partition $\cE=\{E_1,E_2,\dots,E_m\}$ of $E(G)$. We define $\overrightarrow{\Gamma}_{f,\cE}$ to have vertex set $\cE$ and there is a directed edge from $E_i$ to $E_j$ if and only if there exists $e\in E_i$ with $f(e_i)\in E_j$. If all $E_i$s consist of a single edge, then we will drop $\cE$ from the subscript and write $\overrightarrow{\Gamma}_f$. An edge subset $E'\subset E(G)$ is $f$-free if and only if the corresponding vertices form an independent set in $\overrightarrow{\Gamma}_f$.

We now establish the useful Lemma~\ref{superlem} via the following general argument. For any $f\in F_n$ without $f$-fixed $G'\cong G$, the number $m$ of $f$-shifted edges must be at least $\binom{n}{2}-\ex(n,G)$. If there is no $f$-free copy of $H$, then for any $H'\cong H$ in the graph of shifted edges, we must have an $f$-shifted edge $e$ with $e,f(e)$ both belonging to $H'$. We introduce the following two quantities: $C(n,H)$ is the maximum number of subgraphs $H'\cong H$ in $K_n$ that contain a fixed pair of edges (this is a maximum of two values, one for adjacent pairs of edges, one for non-adjacent ones), and $S(n,m,H)$ (the \textit{supersaturation} function) is the minimum number of copies of $H$ in an $n$-vertex graph with $m$ edges. This yields the following statement.

\begin{lemma}\label{superlem}
If $f\in F_n$ admits neither fixed copies of $G$, nor free copies of $H$, then the following inequality holds:
    \begin{equation}\label{super}
    \left(\binom{n}{2}-\ex(n,G)\right)\cdot C(n,H)\ge S(n,\binom{n}{2}-\ex(n,G),H).\tag{*}
\end{equation}
\end{lemma}

Observe that (\ref{super}) always holds if both $G$ and $H$ have chromatic number at least 3, as in this case $2\ex(n, G)$ and $2\ex(n, H)$ exceed $\binom{n}{2}$, and so $\binom{n}{2}-\ex(n,G)\le \ex(n,H)$ and thus $S(n,\binom{n}{2}-\ex(n,G),H)=0$. On the other hand, (\ref{super}) does not hold for large enough $n$ if at least one of $G$ or $H$ is bipartite. This follows from two basic results of extremal graph theory. The Erd\H os--Stone--Simonovits theorem \cite{ES} states that $\ex(n,G)=(1-\frac{1}{\chi(G)-1}+o(1))\binom{n}{2}$. Moreover, Erd\H os and Simonovits showed \cite{ES83} that for any graph $H$ and positive constant $\varepsilon$, there exists a constant $\delta>0$ such that the  inequality $S(n,\ex(n,H)+\varepsilon n^2,H)\ge \delta n^{|V(H)|}$ holds. As $C(n,H)\le n^{|V(H)|-3}$, the right-hand side of (\ref{super}) is $O(n^{|V(H)|-1})$. On the other hand, as at least one of $G$ and $H$ is bipartite, the Erd\H os--Stone--Simonovits theorem implies that $\binom{n}{2}-\ex(n,G)\ge \ex(n,H)+\varepsilon n^2$ for some positive $\varepsilon$. Therefore, the result of \cite{ES83} implies that the right-hand side of (\ref{super}) is of order $n^{|V(H)|}$. So (\ref{super}) cannot hold for all $n$, and thus will yield a bound on $m(G,H)$.

\vskip 0.2truecm

%For recent references to supersaturation, see \cite{CNR, FMS, GNV}.

\vskip 0.2truecm

According to a theorem of Erd\H os \cite{E69}, there exists a constant $\alpha>0$ such that any $n$-vertex graph $G$ with $\ex(n,K_r)+1$ edges has an edge $e$ contained in at least $\alpha n^{r-2}$ copies of $K_r$ in $G$. Applying this repeatedly, we obtain the following supersaturation result.

\begin{theorem}[Erd\H os \cite{E69}]\label{erdos}
Let $r\ge 3$. There exists an $\alpha_r>0$ such that for all positive integers $E$, we have
$S(n,\ex(n,K_r)+E,K_r)\ge \alpha_r\cdot E\cdot n^{r-2}$.
\end{theorem}

For the particular case $r=3$ we shall also use the following more detailed result.

\begin{theorem}[Moon, Moser \cite{MM}]\label{mm}
For any pair $n,m$ of integers, we have $S(n,m,K_3)\ge \frac{4m}{3}(\frac{m}{n}-\frac{n}{4})$. In particular, $S(n,m,K_3)>m$ if $m>\frac{n^2}{4}+\frac{3n}{4}$. Furthermore, if a graph $H$ on $n$ vertices and $m$ edges contains only $\frac{4m}{3}(\frac{m}{n}-\frac{n}{4})$ triangles, then $H$ is regular and every pair of adjacent edges in the complement $H^c$ of $H$ is contained in a triangle of $H^c$. 

Also,  $S(n,h,K_4)\ge \frac{h(4h-n^2)(3h-n^2)}{6n^2}$.  
\end{theorem}

Let us mention that Lov\'asz and Simonovits \cite{LS} determined $S(n,\ex(n,K_r)+t,K_r)$ as long as $t=o(n^2)$ and much later Reiher \cite{R} settled the case $t=\Theta(n^2)$.

\section{The appearance of $f$-fixed copies of $G$ or $f$-free/exclusive copies of $H$: $m(G,H)$ and $m^*(G,H)$}

The theorems in this section supply effective bounds, and sometimes exact values, for $m(G,H)$ and $m^*(G,H)$.

\subsection{\textbf{\boldmath Results about $m(G,K_r)$}}\mbox{}\par\noindent

We start with the proof of Theorem \ref{mtkr}. Let us first note that $C(n,K_r)=\binom{n-3}{r-3}$, in particular $C(n,K_3)=1$.

\begin{proof}[\textbf{Proof of Theorem \ref{mtkr}}]
To see the lower bound, let $n=(r-1)(k-1)$ and consider the subgraph $H=(r-1)K_{k-1}$ of $K_n$. Define the mapping $f\in F_n$ such that $f(e)=e$ if and only if $e\in E(H)$. Then there exists no $f$-fixed copy of $T$, neither any $f$-shifted (and thus no $f$-free) copy of $K_r$. Let us mention that an alternative proof would be to apply Proposition \ref{mzR} to an old result of Chv\'atal \cite{Ch} stating $R(T_k,K_r)=(r-1)(k-1)+1$.

\medskip 

To see the upper bound (2a), let $n:=(r-1)(k-1) + x$ (with an $x$ to be defined later) and consider $f\in F_n$. By the assumption on $\ex(n,T)$, the number of $f$-shifted edges is at least
\[
\binom{n}{2}-\frac{k-2}{2}n=\frac{n(n+1-k)}{2}\ge \frac{[(r-1)(k-1)+x][(r-2)(k-1)+x]}{2}.
\]
By Tur\'an's theorem, we have $$\ex(n,K_r)\le \left(1-\frac{1}{r-1}\right)\binom{n}{2}\le \frac{[(r-1)(k-1)+x][(r-2)(k-1)+\frac{(r-2)x}{r-1}]}{2}.$$ So the number of $f$-shifted edges is at least $\ex(n,K_r)+\frac{x}{2(r-1)}n$. Theorem~\ref{erdos} and Lemma~\ref{superlem} (\ref{super}) imply  $$\binom{n}{2}\binom{n-3}{r-3}\ge\frac{x\alpha}{2(r-1)}n^{r-1},$$  which does not hold if $x$ is a large enough constant depending only on $r$.

\medskip 

The upper bound of (2b) follows from Lemma~\ref{superlem} (\ref{super}), as for any $f\in F_{2k+2}$, if there do not exist $f$-fixed copies of $T$, the left-hand side is at most the number $m$ of $f$-shifted edges, where $m\ge \binom{2k+2}{2}-\ex(n,T)\ge (k+1)(k+3)$. So with $m\ge (k+1)(k+3)$ and $n=2k+2$, we have $\frac{m}{n}-\frac{n}{4}\ge \frac{(k+1)(k+3)}{2k+2}-\frac{2k+2}{2}=1$. Therefore Theorem \ref{mm} implies that the right-hand side is at least $\frac{4m}{3}>m$ and so Lemma~\ref{superlem} (\ref{super}) does not hold, implying $m(T,K_3)\le 2k+2$ for any EST. Doing the same computation for $n=2k+1$, the left-hand side of Lemma~\ref{superlem} (\ref{super}) is still $m$. The number of $f$-sifted edges is at least $\binom{n}{2}-\frac{k-2}{2}n=\frac{1}{2}n(n+1-k)$. By Theorem \ref{mm} the lower bound $S(n,m,K_3)$ is at least
\[
\frac{4m}{3}\left(\frac{m}{n}-\frac{n}{4}\right)=\frac{4m}{3}\left(\frac{n+1-k}{2}-\frac{n}{4}\right)=\frac{4m}{3}\left(\frac{k+2}{2}-\frac{2k+1}{4}\right)=\frac{4m}{3}\cdot \frac{3}{4}=m.
\] The furthermore part of Theorem \ref{mm} implies that either the left hand side is strictly more than $m$ and we got our contradiction and $m(T,K_3)\le 2k+1$ or the complement of the graph of shifted edges, i.e. the graph of fixed edges must be the union of cliques. In such a $T$-free graph with $\frac{k-2}{2}n$ edges, the cliques must have size ${k-1}$. In particular, $k-1$ should divide $2k+1$, so either $k=2$ and $T=K_2$ or $k=4$ and $T=K_{1,3}$ or $T=P_4$.

\medskip 

The lower bound of (2c) was obtained in \cite{AC} for $K_2$ and $K_{1,3}$. Here we prove it for $P_4$ and $K_{1,3}$. For the proof, it is enough to define $f \in F_9$ with no $f$-fixed copies of $P_4$ and $K_{1,3}$, neither $f$-free copies of $K_3$. Consider a partition $A\cup B\cup C=V(K_9)$ with $|A|=|B|=|C|=3$. For any edge $e$ contained in any of the parts, we let $f(e)=e$; and for all crossing edges, we will have $f(e)\neq e$. This ensures that there is no $f$-fixed copy of $K_{1,3}$ or $P_4$. There are 27 crossing edges, and there are 27 crossing triangles (a crossing triangle is a triangle that contains one vertex from each of $A$, $B$, and $C$), and using Hall's theorem it is easy to see that there is a one-to-one mapping $\eta$ between crossing edges and crossing triangles in such a way that $e \subset \eta(e)$ for all crossing edges $e$. So if we let $f(e)$ be one of the other two edges of $\eta(e)$, then there do not exist any $f$-free copies of $K_3$.

\medskip

The lower bound of (2d) is a
special case of the lower bound of part (1) of Theorem 1. To obtain the upper bound, let $n$ be an integer and let $f \in F_n$  be a function such that $K_n$ contains neither $f$-fixed copies of $T$, nor $f$-free copies of $K_4$.

Since $K_n$ contains no $f$-fixed copy of the Erd\H os--S\'os tree $T$, the number of $f$-fixed edges is at most $(k- 2)n/2$, hence the number $h$  of $f$-shifted edges satisfies $h\ge  \binom{n}{2}  - \frac{(k-2)n}{2}$.  Let $H$ be the subgraph of $K_n$ consisting of the $f$-shifted edges. Let $S$ be the set of all ordered pairs $( e, K)$, where $K$ is a $K_4$  in $H$, and $e ,f( e)\in K$. 
By Theorem \ref{mm} for $K_4$, a graph with $n$ vertices and $h$ edges contains at least  $\frac{h(4h-n^2)(3h-n^2)}{6n^2}$ copies of $K_4$.

Since $H$ contains no $f$-free $K_4$, all copies of $K_4$ in $H$ appear at least once as a second coordinate of a pair in $S$, hence $|S|\ge   \frac{h(4h-n^2)(3h-n^2)}{6n^2}$.
Clearly every edge $e \in E(H)$ 
can appear in at most $n-3$ pairs of $S$, thus $h(n-3)\ge|S|$. 
Combining the last two inequalities we conclude that $(n-3)h \ge \frac{h(4h-n^2)(3h-n^2)}{6n^2}$  or equivalently   $12h^2 - 7hn^2 +n^4- 6n^3 + 18n^2 \le  0$. Solving for $h$ we get $h \le \frac{7n^2 +\sqrt{n^4+288n^3-864n^2}}{24}  <  \frac{7n^2 +  n^2 +144n}{24} = \frac{n^2}{3}  + 6n$.

Combining the above upper and lower bounds on $h$,  we conclude that  $n^2\!/3 +6n \ge h\ge \binom{n}{2}  -\frac{(k-2)n}{2}$,  and $n \le 3k +33$.
\end{proof}

The following theorem concerning $m(T,2K_2)$ is the complement of a theorem concerning $h(n,2K_2)$ proved in \cite{CPTV}.

\begin{thm}
For any Erd\H os--S\'os tree $T$ on $k \ge 7$ vertices we have $m(T,2K_2) = k+1$.
\end{thm}

\begin{proof}
To prove the lower bound consider the following $f\in F_k$: fix $x\in V(K_k)$, let $f(e)=e$ if $x\notin e$, and for all $y$ define $f(xy)=(xy')$ with some $y'\neq y$. As all shifted edges form a star, $f$ does not admit free copies of $2K_2$ and there is no $f$-fixed copy of $T$ on $k-1$ vertices. This shows $m(T,2K_2) \ge k + 1$.

To prove the upper bound for $k+1 \ge 8$, consider $f\in F_{k+1}$. 
If the $f$-shifted graph has at most $k$ edges, then the number of $f$-fixed edges is at least $(k +1)k/2  - k  =  k(k-1)/2 > (k+1)(k-2)/2\ge \ex(k+1,T)$. So $f$ admits a fixed copy of $T$.

In \cite{CPTV}, the authors proved that  $h(n,2K_2) = n$ holds for any $n \ge 7$, and the only extremal graph is $K_{1,n-1}$ plus an edge. Therefore, if there are more than $k +1$ $f$-shifted edges, then $f$ admits a free copy of $2K_2$. Furthermore, if for some $f\in F_{k+1}$ there are exactly $k+1$ $f$-shifted edges without free copies of $2K_2$, then the shifted edges must form a $K_{1,k}$ plus an edge. Hence the $f$-fixed graph is $K_k$ minus an edge,  and this contains copies of all trees on $k$ vertices.
\end{proof}

The next theorem uses Lemma~\ref{superlem} to derive an upper bound on $m(G,K_r)$ in case $G$ is bipartite and the exact value or a good upper bound on $\ex(n,G)$ is known.

\begin{thm}\label{mgkr}
\
\begin{enumerate}
    \item
    Let $G$ be a bipartite graph, and assume that $n$ is an integer satisfying $\lceil \frac{n^2}{4} \rceil-\frac{5n}{4} > \ex(n,G)$.  Then  $m(G,K_3) \le n$ holds.
    
\medskip 
    
    \item
    For any $r\ge 3$, there exists a constant $C=C_r$ such that if $G$ is a bipartite graph  with $\ex(n,G)<\frac{1}{r-1}\binom{n}{2}-C\cdot n$, then $m(G,K_r)\le n$ holds. 
\end{enumerate}

\end{thm}

\begin{proof}
To see (1), let $n$ be an integer with $\lceil \frac{n^2}{4}\rceil-\frac{5n}{4}>\ex(n,G)$ and let $f\in F_n$. Suppose $K_n$  contains neither $f$-fixed copies of $G$ nor $f$-free triangles. As $\lceil \frac{n^2}{4} \rceil-\frac{5n}{4}> \ex(n,G)$, we have $\binom{n}{2}-\ex(n,G)>\lfloor \frac{n^2}{4}\rfloor+\frac{3n}{4}$. 
Plugging into Lemma~\ref{superlem} (\ref{super}) and using Theorem \ref{mm} stating $m>S(n,m,K_3)$ if $m>\frac{n^2}{4}+\frac{3n}{4}$ and $C(n,K_3)=1$, we obtain
a contradiction if $\binom{n}{2}-\ex(n,G)>\lfloor \frac{n^2}{4}\rfloor+\frac{3n}{4}$ is equivalent to $\binom{n}{2}-\ex(n,G)> \frac{n^2}{4}+\frac{3n}{4}$. For even $n$ the two left-hand sides are the same, while if $n$ is odd, then $\lfloor \frac{n^2}{4}\rfloor+\frac{3n}{4}$ is either an integer or an integer plus $1/4$ and always $1/4$ smaller than $\frac{n^2}{4}+\frac{3n}{4}$. So any integer, in particular $\binom{n}{2}-\ex(n,G)$ that is strictly larger than  $\lfloor \frac{n^2}{4}\rfloor+\frac{3n}{4}$ is also strictly larger than  $\frac{n^2}{4}+\frac{3n}{4}$. This completes the proof of (1).

Similarly, for the proof of (2) fix $n$ large enough and $f\in F_n$. Note that $C(n,K_r)=\binom{n-3}{r-3}$, so the right-hand side of Lemma~\ref{superlem} (\ref{super}) is at most $\binom{n}{2} \binom{n-3}{r-3}$. On the other hand, using Tur\'an's theorem on $\ex(n,K_r)$, the bound on $\ex(n,G)$ yields $\binom{n}{2}-\ex(n,G)>(1-\frac{1}{r-1})\binom{n}{2}+C\cdot n\ge \ex(n,K_r)+C\cdot n$. Theorem \ref{erdos} implies that the right-hand side of Lemma~\ref{superlem} (\ref{super}) is at least $C\cdot \alpha \cdot n^{r-1}$ for the $\alpha$ given there. If $C$ is large enough, this is a contradiction. 
\end{proof}

%\begin{thm}
%  Let $T$  be an EST with $k$ vertices for some $k \ge 2$. Then 

%\begin{enumerate}[a.]
%    \item 
%    $m(T ,K_3) \le  2k + 2$ and  if $T\neq  K_2, K_{1,3} ,  P_4$,  then  $m(T,K_3) \le 2k +1$.
%    \item
%     $m(T,K_3)  \ge 2k -1$.
%     \item
%     if $T \in  \{K_2,K_{1,3}, P_4 \}$,  then  $m(T, K_3)  =  2k +2$.
%\end{enumerate}  
%\end{thm}

\subsection{\textbf{\boldmath Results about $m(G,K_{1,r})$}}\mbox{}\par\noindent

We need the next two lemmas to prove Theorem~\ref{th:m-tree}.

\begin{lemma}\label{euler}
    Suppose $K_n$ can be partitioned into two graphs $H_1$, $H_2$ such that $H_1$ is $G$-free, and $\Delta(H_2)\leq 2r-2$. Then $m(G,K_{1,r})>n$ holds. 
\end{lemma}

\begin{proof}
    We need to define an $f\in F_n$ that admits no fixed copy of $G$ and no free copy of $K_{1,r}$. We let $f(e)=e$ if and only if $e\in E(H_1)$ and thus by the assumption on $H_1$, $f$ does not admit a fixed copy of $G$. The number of vertices that have odd degrees in $H_2$ is even. Therefore if we add a new vertex $w$ to $H_2$ and connect $w$ to all odd-degree vertices, the obtained graph $H'_2$ will only have even degrees and thus contains an Eulerian walk
    (or each of its components does if all degrees in $H_2$ are even and no $w$ has been added).

    Let $e_1,e_2,\dots,e_{h}$ be the edges in the order of the walk. Then for $e_i\in E(H_2)$, we define $f(e_i)=e_{i+1}$ with addition considered modulo $h$ if $e_{i+1}\in E(H_2)$; otherwise $w$ is an endpoint
of $e_{i+1}$ and we pick $e^*\neq e_i$ arbitrarily and let $f(e_i)=e^*$. Then at every vertex $x$, there are some $a\le r-1$ incoming and the same number of outgoing edges of the walk. By definition of $f$, an $f$-free star cannot contain both of a consecutive pair of edges (either because it would ruin the $f$-free property, or because one of these edges is connected to $w$ so not in $E(H_2)$). Therefore, the largest $f$-free star contains only $r-1$ edges.
\end{proof}

The next lemma is the well-known Frobenius--Sylvester theorem (see \cite{KZ} for a recent reference).

\begin{lemma}\label{lincomb}
    If $a$ and $b$ are coprime, then any integer larger than $ab-a-b$ can be expressed as $xa+yb$ with $x,y$ being non-negative integers.
\end{lemma}

\begin{proof}[\textbf{Proof of Theorem \ref{th:m-tree}}]
    The upper bound of (1) follows from Proposition \ref{upperk1r} as for an EST $T$ we have $\binom{n}{2}-\ex(n,T)\ge \binom{n}{2}-\frac{k-2}{2}n=\frac{n}{2}(n-k+1)>n(r-1)$ for $n=k+2r-2$. The upper bound of (2) follows similarly, as for a strong EST $T$, if $k-1$ does not divide $2(r-1)$ and so it does not divide $k+2r-3$. Therefore, for $n=k+2r-3$ we have $\binom{n}{2}-\ex(n,T)>\binom{n}{2}-\frac{k-2}{2}n=n(r-1)$.
    
    To prove the lower bound of (1), suppose $(k -1) \mid 2(r-1)$  and apply Lemma \ref{euler} with $n=k -1+ 2(r-1)$ as follows.   Split $V(K_{k+2r-3})$  into  $t =  1+ \frac{2r-2}{k-1}$  subsets $ V_1,\dots,V_t$ of cardinality $k-1$ each, and let $H_1$ be the disjoint union of cliques on the $V_i$s and $H_2$ be the graph on the remaining edges. As the components of $H_1$ are of size less than $k$, $H_1$ is $T$-free, and $H_2$ is $2(r-1)$-regular, so $\Delta(H)\le 2r-2$.

%\medskip 

 %   To prove (2) suppose $n =  k +2r  - 3$. Since $2(r-1) \not =  0$ (mod $k-1$),  it follows that $n \not = 0$ (mod $k-1$).   If there is no $f$-fixed copy of $T$, then as $T$ is a strong EST,  the number of $f$-shifted edges is strictly larger than $\frac{n(n-1)}{2}  - \frac{(k-2)n}{2}  = n(r-1) $.  Hence by the definition of $d_{sh}^{*}(H)$ there is a vertex with $d_{sh}(x) \ge r$ and an $f$-free $K_{1,r}$ occurs.

\medskip 

    To prove the lower bound of (3), we apply Lemma \ref{euler} with $n=k -3+ 2(r-1)$ as follows. We set $a=k-1$ and $b=k-3$, and observe that their greatest common divisor is either 1 or 2. In both cases, by Lemma \ref{lincomb}, any even number at least $k^2$ can be expressed as a non-negative linear combination of $a$ and $b$. In particular, $2(r-1)=x(k-1)+y(k-3)$ and thus $n=x(k-1)+(y+1)(k-3)$. We split $K_n$ into $H_1$, $H_2$ with $H_1$ being $xK_{k-1}+(y+1)K_{k-3}$ and $H_2=K_n\setminus H_1$. Then $H_1$ is $T$-free as all its components have size less than $k$, and the degrees in $H_2$ are $2(r-1)$ and $2(r-1)-2$.

\medskip 

    Finally, to prove the lower bound of (4), we apply Lemma \ref{euler} with $n=k -2+ 2(r-1)$. This time we set $a=k-1$ and $b=k-2$, which are clearly coprime. So if $r\ge k^2$, then $2(r-1)=xa+yb$ and thus $n=xa+(y+1)b$. Let $K_{k-1}^-$ be the graph we obtain from $K_{k-1}$ with a 1-factor removed (here we use the assumption that $k$ is odd, and thus $k-1$ is even). We split $K_n$ to $H_1,H_2$ with $H_1=xK_{k-1}^-+(y+1)K_{k-2}$ and $H_2=K_n\setminus H_1$. Then $H_1$ is $T$-free as above and $H_2$ is $2(r-1)$-regular.
\end{proof}

We need the next two lemmas to prove Theorem~\ref{*tk}.

\begin{lemma}\label{excluk1r}
\begin{enumerate}[(i)]
    \item Let $G$ be a graph in which no edge is incident to all other edges. Suppose $H$ contains $G$, and let $f\in F_{H,0}$. If $\Delta(G) \ge 5r - 4$, then there exists in $G$ an $f$-exclusive copy of $K_{1,r}$.

    \item Similarly, if $G$ is a graph containing $M = (5t-4)K_2$ and $f : E(G)\rightarrow E(G)$ is strong-shifted on $M$, then $G$ contains an $f$-exclusive copy of $tK_2$.
\end{enumerate}
\end{lemma}

\begin{proof}
    Let $v$ be a vertex of $G$ with $d(v) \ge 5r - 4$ and let $e_1, e_2, \dots , e_t$ be the edges incident to $v$ in $G$. We will apply Lemma \ref{outdeg} to the directed graph $\overrightarrow{D}$ whose vertex set is $\{e_1,e_2,\dots,e_t\}$, with an arc directed from $e_i$ to $e_j$ if and only if  $|f(e_i)\cap   e_j | \ge 1$. 
As $f\in F_{H,0}$, it follows that  $v\notin f(e _j)$ for all $j = 1,\dots ,t$. Hence the maximum out-degree in $\overrightarrow{D}$ is at most 2.

Lemma \ref{outdeg}  implies the existence of an independent set of cardinality at least $\lceil \frac{d(v)}{5}\rceil \ge \lceil \frac{(5r - 4)}{5} \rceil = r$ and these edges form an $f$-exclusive copy of $K_{1,r}$.  

The proof of the second statement is basically the same.
\end{proof}

\begin{lemma}\label{*k1r}
    If  $r \ge 2$  and  $\binom{n}{ 2} -\ex(n,G)  >  (5r-5)n/2$, then  $m^*(G,K_{1,r}) \le  n$.
Further, if $\binom{n}{2}  - \ex(n,G) >  \ex(n, (5t - 4)K_2 )$, then $m^*(G,tK_2)\le n$.
\end{lemma}

\begin{proof}
    Let $f\in F_n$ be such that $|f(e)\cap e|$ is either 0 or 2 for all edges $e$. If $f$ does not admit a fixed copy of $G$, then the graph of strong-shifted edges has more than $(5r-5)n/2$ edges and hence has maximum degree at least $5r-4$, and Lemma~\ref{excluk1r}(i) applies. This proves the first statement of Lemma~\ref{*k1r}. The proof of the second statement is similar.
\end{proof}

Now we are in a position to prove Theorem \ref{*tk}.

\begin{proof}[Proof of Theorem \ref{*tk}]
    For the upper bound of (1), observe that for an EST on $k$ vertices, the condition $\binom{n}{2}-\ex(n,T)\ge\binom{n}{2} - \frac{(k-2)n}{2}  > \frac{(5r -5 )n}{2}$ of Lemma \ref{*k1r} simplifies to $n>5r+k-6$.

    To obtain the lower bound of (2), we need to construct $f\in F_n$ with $n=2r+k-3=2(r-1)+k-1$ satisfying $|e\cap f(e)|\neq 1$ for all $e$ such that $f$ does not admit any fixed copies of $T$, nor any exclusive copies of $K_{1,r}$. By the divisibility condition, $n=s(k-1)$ holds for some $s$. Consider a copy $K$ of $sK_{k-1}\subset K_n$ and define $f(e)=e$ for all $e\in E(K)$. Clearly, no fixed copies of $T$ are created and the edges of $E_n\setminus E(K)$ form an $(n-k+1)$-regular graph $G$. 

    Since $k-1$ is even, the edge set between any two components of $K$ admits a cycle decomposition because it is isomorphic to $K_{k-1,k-1}$.
    Note that each vertex is incident with exactly $r-1$ of the cycles obtained in these decompositions. In this way, $E(G)$ is decomposed into a collection of even cycles, say $C^1, C^2,\dots$.

    If $E(C^j)$ is $e_1,e_2,\dots,e_{2t}$, then we let $f(e_i)=e_{i+2}$ where addition is modulo $2t$. By this definition, for every $e\in E_n\setminus E(K)$ we have $e\cap f(e)=\emptyset$. Also, if $S$ is a copy of $K_{1,r}$ in $G$ with center $v$, then there exist two edges $e,e'$ of $S$ that belong to the same 
    cycle $C^i$ incident with $v$, 
    and thus $e=e_h$ and $e'=e_{h+1}$ are consecutive in some cycle. But then $|f(e)\cap e'|=1$, so $S$ is not an $f$-exclusive copy of $K_{1,r}$. 
\end{proof}

The next lemma and proposition are from \cite{AC}; we restate and prove them for the sake of being self-contained. Proposition~\ref{upperk1r} will be used below, in particular in Theorems~\ref{properties} and~\ref{th:exact-families}.
For $f\in F_n$ and $x\in V(K_n)$, let the \textit{shifted degree of $x$} be $d_{sh}(x):=|\{e: x\in e, \ x \notin f(e) \}|$, and let us define \textit{the average shifted degree of $f$} to be $d_{sh}^*(f):=\frac{1}{n} \sum_{x}d_{sh}(x)$.

\begin{lemma}\label{s*}
    If $d_{sh}^*(f)>r-1$, then $f$ admits an $f$-free $K_{1,r}$.
\end{lemma}

\begin{proof}
    There exists an $x$ with $d_{sh}(x)>r-1$ and thus $d_{sh}(x)\ge r$. The corresponding edges incident with $x$ form an $f$-free $K_{1,r}$.
\end{proof}

\begin{proposition}\label{upperk1r}
    For a graph $G$ and an integer $r$, if $n$ is an integer such that $\binom{n}{2}  -  \ex(n,G)  >  n(r-1)$, then $m(G,K_{1,r}) \le n$.
\end{proposition}

\begin{proof}
    Let $f\in F_n$. By assumption, the graph $H$ of non-fixed edges has $e(H)>n(r-1)$.  Then as each $xy\in E(H)$ contributes 1 to at least one of $d_{sh}(x)$ and $d_{sh}(y)$, we have $e(H)\le \sum_{x}d_{sh}(x)$, and so $(r-1)<d_{sh}^*(f)$. Lemma \ref{s*} completes the proof.
\end{proof}

In the following theorems and propositions, we develop, using the notion of the Deck of a graph, a method that allows us to exactly determine $m(G,K_{1,r})$ for many graphs $G$ with chromatic number at least 3 and widely expand the range of results proved in \cite{AC}.

\medskip

Given a graph $G$, let the \textit{Deck of $G$} be defined as $D(G):=\{ G \setminus \{v\}: v\in V(G) \}$.

\noindent 
\begin{thm}\label{properties}

\ 

\begin{enumerate}
    \item \label{gener}
    For any graphs $G$  and $Q$, we have $m(G, Q) \le \min\{m(H, Q) : H \in D(G)\} + \lfloor  \frac{2h(n,Q)}{n}\rfloor  +1$.

    \item \label{deck} For any graph $G$ and $r \ge 1$, we have $m(G,K_{1,r})  \le \min \{ m(H,K_{1,r}): H \in D(G) \} +  2r-1.$

\smallskip 

    \item \label{chi} If  $\chi(G)   = k \ge 3$ and  $r \ge 1$,  then $m(G,K_{1,r})  \ge (k-1)(2r-1) +1$.

\smallskip 

    \item \label{join} Let $\chi(G) = k \ge 3$ and suppose $m(G,K_{1,r})  = (k-1)(2r-1) +1$.  Then for $t \ge 1 $, we have $m(G+ K_t,K_{1,r})  =  (k+t -1)(2r-1) +1$.

     \end{enumerate} 
 \end{thm}

 \begin{proof}
To prove (\ref{gener}), write $m :=  \min\{m(H, Q) : H \in D(G)\}$. Let $f\in F_n$,  where $n = m + \lfloor\frac{2h(n,Q )}{n}\rfloor +1$.
If there is a vertex $v$ with at least $m$ fixed edges incident with it, then in the restriction of $f$ to $N[v]$ we have either a fixed copy of $G$ or a free copy of $Q$.
So we may assume that every vertex has at most $m - 1$ fixed incident edges, hence at least $n - 1 - (m - 1) = n - m =  \lfloor\frac{ 2h(n, Q)}{n }\rfloor +1$ non-fixed edges. Then the total number of shifted edges is at least $\lfloor\frac{ 2h(n, Q)+n}{2}\rfloor > h(n, Q)$ and so by definition there must exist an $f$-free copy of $Q$. 
 
 To prove (\ref{deck}), write $m : = \min \{ m(H,K_{1,r}): H \in D(G) \}$.  Let $f \in F_n$,  where $n  =  m +2r  - 1$. If there is a vertex  $v$  with at least  $m$ fixed edges incident with it, then together with $N(v)$ we have either a fixed copy of $G$ or a free $K_{1,r}$.
So we may assume that every vertex has at most  $m -1$ fixed incident edges, hence at least $n -1 - (m-1)  =  n-m  = 2r-1$ non-fixed edges. Every non-fixed edge contributes at least 1 to $\sum_xd_{sh}(x)$, so we obtain $d_{sh}^*(f)  \ge   \frac{1}{n}\frac{(2r-1)n}{2} > r  - 1$ and therefore $f$ admits a free $K_{1,r}$ by Lemma \ref{s*}.

The claim of (\ref{chi})  is trivial for $r = 1$. Let $r\ge 2$,  $n = (k-1)(2r-1)$, and consider $k -1$ copies of $K_{2r -1}$. On each copy of $K_{2r -1}$ define $f =  f^*$ that realizes the lower bound for $g(K_{1,r}) = 2r$ (see Theorem \ref{hkr1} in Section 5), what we know from \cite{AC}.  All the edges not belonging to these copies of $K_{2r -1}$ are $f$-fixed. By construction, $f$ does not admit free copies of $K_{1,r}$ and fixed copies of $G$ as the subgraph determined by the $f$-fixed edges is $(k-1)$-chromatic.

To see the lower bound in (\ref{join}), observe  $\chi(G+ K_t) = \chi(G) +t =  k+t$. Hence by Theorem \ref{properties} (\ref{chi}) we obtain $m(G+ K_t, K_{1,r}) \ge (k +t - 1)(2r-1) +1$. The upper bound follows by induction on $t$, applying Theorem \ref{properties} (\ref{deck}).
 \end{proof}

Along the lines of the proof of Theorem \ref{properties} (\ref{gener}), one can show the following analogous statement.

\begin{proposition}

 For any graphs $G$  and $Q$, we have $m^*(G, Q) \le \min\{m^*(H, Q) : H \in D(G)\} + \lfloor \frac{ 2s(n,Q)}{n}\rfloor  +1$.
\end{proposition}

%\begin{proposition}
%    Let $\chi(G) = k \ge 3$ and suppose $m(G,K_{1,r})  = (k-1)(2r-1) +1$.  Then for $t \ge 1 $, we have $m(G+ K_t,K_{1,r})  =  (k+t -1)(2r-1) +1$.
%\end{proposition}

%\begin{proof}

%\end{proof}

Let us introduce the notation $K_k^{-t}$ for the graph obtained from $K_k$ by removing the edges of a $t$-clique. In the next two theorems, we determine $m(G,K_{1,r})$ for members $G$ of several graph classes. The lower bound will be obtained by Theorem \ref{properties} and we will apply Proposition \ref{upperk1r}. To this end, we will need the exact value of $\ex(n,G)$. This is known for all color-critical graphs by a result of Simonovits \cite{S}. But his result is only proved for large enough $n$, while we will need the parameter for small $n$, so instead of using this general theorem, we will use specific results. 

\begin{thm}\label{th:exact-families}
\

    \begin{enumerate}
    \item
    For $k \ge t+2$, $r \ge \frac{3t}{2}-1$, we have  $m(K_k^{-t}, K_{1,r})  = (k-t)(2r -1) +1$.  
    \item
    For $k \ge 3$, $r\ge  k -1$, $t \ge 0$, we have  $m(C_{2k-1}+ K_t, K_{1,r}) = (t + 2)(2r - 1) +1$.
    \end{enumerate}
\end{thm}

\begin{proof}
    Let us start with the proof of (1). First observe that it is enough to prove the statement for $k=t+2$, as for larger $k$ it then follows from $K_k^{-t}=K_{t+2}^{-t}+ K_{k-t+2}$ and Theorem \ref{properties}~(\ref{join}).

    So it is enough to prove that for  $r\ge \frac{3t}{2}-1$, we have $m(K_{t+2}^{-t}, K_{1,r})  = 2(2r -1)  +1 = 4r - 1$. The lower bound  follows from Theorem \ref{properties} (\ref{chi}) as $\chi(K_{t+2}^{-t}) = 3$. Also, by a result of Edwards \cite{Ed} and Had\v ziivanov and Nikiforov \cite{HN}, we have $\ex(n,K_{t+2}^{-t})=\lfloor \frac{n^2}{4}\rfloor$ if $n\ge 6t$. To see the upper bound, let $f\in F_{n}$ with $n=4r-1$ and apply Proposition \ref{upperk1r}:
    %and observe that either the number of $f$-shifted edges is at least $n(r-1)+1$, and then by Proposition \ref{hkr1} (1), we have an $f$-free copy of $K_{1,r}$ or the number of $f$-fixed edges is at least ű
    $$\frac{n(n-1)}{2}   - n(r-1)=  \frac{n(n +1 - 2r)}{2}  = (4r-1)r > \frac{n^2}{4} =  \ex(n,K_{t+2}^{-t}).$$ So $m(K^{-t}_{t+2},K_{1,r})\le n$.

    As in (1), for (2) it is enough to consider $t=0$ as the $t>0$ case follows by Theorem \ref{properties} (\ref{join}). The lower bound again follows from Theorem \ref{properties} (\ref{chi}) as $\chi(C_{2k-1})=3$. Also, $\ex(n, C_{2k-1})=\lfloor \frac{n^2}{4}\rfloor$ if $n\ge 4k-4$ (see \cite{FG}), and the same computation as for the upper bound of (1) and Proposition \ref{upperk1r} completes the proof.
\end{proof}

\begin{thm}\label{th:exact-subgraphs}
\begin{enumerate}
    \item
     Let $t \ge 2$, $r \ge  3t/2 -1$, and let $G$ be a graph on $k \ge t +2$  vertices with $\chi(G) = k - t +1$ such that $G$ is  contained in $K_k^{-t}$.   Then  $m(G,K_{1,r}) =(k-t)(2r -1) +1$.   
\item
Let $k\ge 3$, $r\ge  k-1$, $t \ge 0$, and let $G$ be a graph on $t + 2k  - 1$ vertices with $\chi(G) = t +3$ such that  $G$ is  contained in $C_{2k-1} + K_t$.   Then  $m(G,K_{1,r}) =(t +2)(2r -1) +1$.  
\item
  There exist positive constants $c_1$ and $c_2$ such that  $2r + c_1\sqrt{r}\le m(C_4  ,K_{1,r}) \le  2r + c_2\sqrt{r}$. 
  \end{enumerate}
  \end{thm}
  
\begin{proof}
    To see (1), as  $\chi(G) =k-t +1 \ge 3$, we can apply Theorem \ref{properties} (\ref{chi}) and Theorem \ref{properties} (\ref{join}) to obtain $$(k-t)(2r-1) +1 \le m(G, K_{1,r}) \le m(K_k^{-t}, K_{1,r})=(k-t)(2r-1) +1.$$

    To see (2), as  $\chi(G) =t +3 \ge 3$, we can apply Theorem \ref{properties} (\ref{chi}) and Theorem \ref{properties} (\ref{join}) to obtain $$(t+2)(2r-1) +1 \le m(G, K_{1,r}) \le m(C_{2k-1} + K_t, K_{1,r})=(t+2)(2r-1) +1.$$

    To see the lower bound of (3) consider a regular $n$-vertex $C_4$-free graph $H$ that has the largest number of edges among all such graphs. By a result of Tait and Timmons \cite{TT}, we know that $e(H)=(\frac{1}{2\sqrt{6}}-o(1))n^{3/2}$ and so all degrees in $H$ are $d=d_n=(\frac{1}{\sqrt{6}}-o(1))\sqrt{n}$. As long as $n-1\le d+2r-2$, we can partition $K_n$ to $H=H_1$ and its complement $H_2$ such the conditions of Lemma \ref{euler} are satisfied and we obtain $n<m(C_4,K_{1,r})$. 
%So for large enough $r$, the constant $c_1$ can be chosen to be $\frac{1}{\sqrt{6}}-o(1)$.
So a suitable positive constant $c_1<\frac{1}{\sqrt{6}}$ can be chosen for large enough $r$.

    For the upper bound of (3), we use the well-known result \cite{ERS} $\ex(n,C_4)=(\frac{1}{2}+o(1))n^{3/2}$. Applying Proposition \ref{upperk1r}, we obtain that a constant $c_2>\frac{1}{\sqrt{2}}$ can be chosen if $r$ is large enough.
\end{proof}

\noindent\textbf{Remark.} Based on extremal results on intersecting odd cycles \cite{kinai} and intersecting cliques $K_1+tK_r$ \cite{chenetal} (and \cite{erdosetal} for $r=2$), similar results on $m(G,tK_{1,r})$ and $m(G,tK_2)$ can be obtained for these graphs $G$.

\medskip

\section{The appearance of $f$-free copies of $G$:  $g(G)$  }

In this section we prove Theorems~\ref{degreethm} and~\ref{gGkm} and additionally we exactly determine $g(tK_2,0)$ and $g(tK_2,1)$.

\begin{proof}[\textbf{Proof of Theorem \ref{degreethm}}]
Note first that $G$ contains $P:=\sum_{i=1}^k \binom{d_i}{2}$
 copies of $P_3$ and $M:=\binom{m}{2}-P$ copies of $2K_2$.

Consider any $f\in F_{n,1}$.
We now view both $G$ and $K_n$ as vertex-labeled graphs.
Then the number of (labeled) copies of $G$ in $K_n$ is
 precisely $\binom{n}{k} k!$.
A copy can be made \textit{not} $f$-free in two ways:
 \begin{itemize}
  \item[$(a)$] there is a $P_3=xyz$ in $G$ with $f(xy)=yz$;
  \item[$(b)$] there is a $2K_2$ in $G$ with two edges $wx$ and $yz$
   such that $f(wx)=yz$.
 \end{itemize}

In case of $(a)$ we need an edge $e\in E_n$ with $e\cap f(e)\neq \emptyset$,
 and once such an $e$ has been chosen, the position of $y$ is determined.
There are two ways to map $\{x,z\}$ onto $(e\cup f(e)) \setminus (e\cap f(e))$,
 and the remaining $k-3$ vertices can be mapped into $V(K_n)$ in
 $\binom{n-3}{k-3}(k-3)!$ different ways.
We have at most $\binom{n}{2}$ choices for $e$ and exactly $P$
 choices for $P_3$, hence the number of non-$f$-free labeled copies
 of $G$ caused by a $P_3$ is at most $\frac{P\cdot n!}{(n-2)\cdot(n-k)!}$
 that we write in the form $A\binom{n}{k} k!$.

On the other hand, for $(b)$ we need an edge $e\in E_n$ with $e\cap f(e)=\emptyset$.
Having chosen such an $e$, for each $2K_2$ of $G$ we have a choice
 between $wx$ and $yz$ to map it onto $e$, and each of them
 can be mapped in two ways onto $e$ and $f(e)$.
Hence $8M$ such mappings exist for a fixed $e$.
Each of them can be extended to a mapping of the entire $G$ in
 $\binom{n-4}{k-4}(k-4)!$ different ways, and the number of choices
 for $e$ is again at most $\binom{n}{2}$.
Consequently, the number of non-$f$-free labeled copies
 of $G$ caused by a $2K_2$ is at most
  $\frac{4M\cdot n!}{(n-2)(n-3)\cdot(n-k)!}$
 that we write in the form $B\binom{n}{k} k!$.

As $K_n$ contains just $\binom{n}{k} k!$ copies of $G$,
 an $f$-free copy of $G$ must occur whenever $A+B<1$ holds.
This inequality means
 \begin{equation}   \label{eq:A+B}
   1 > A+B = \frac{P}{(n-2)} + \frac{4M}{(n-2)(n-3)}
    = \frac{(n-3)\cdot P + 4(\binom{m}{2}-P)}{(n-2)(n-3)} \tag{**}
 \end{equation}
 and the theorem follows.
\end{proof}

As a side-product, from inequality (\ref{eq:A+B}) we also obtain the following estimate.

\begin{proposition}
Let $G$ be a graph with $k$ vertices, $m$
edges, and degree sequence $d_1, d_2, \dots , d_k$. With the notation $P=\sum_{i=1}^k\binom{d_i}{2}$ and $M=\binom{m}{2}-P$ introduced above, we have
 $$
 g(G) \leq
  \max \lbrace |P| , \left\lceil 2\sqrt{|M|} \right\rceil \rbrace
     + 3.
 $$
\end{proposition}

Hence, if $P>0$ and $M>0$ (i.e., $G$ is neither a matching nor a star),
upper bounds more sensitive than the one in Theorem \ref{degreethm}
can be given.

\begin{proof}[\textbf{Proof of Theorem \ref{gGkm}}]
Again, we view $G$ and $K_n$ as vertex-labeled graphs.
Consider any mapping $f\in F_{n,0}$.
As above, the number of (labeled) copies of $G$ in $K_n$ is
 precisely $\binom{n}{k} k!$.
If we fix an $e\in E_n$, an edge $e'\in E(G)$ can be mapped onto $e$ in two ways, and to have a non-$f$-exclusive copy of $G$ we need that a vertex from $V(G)\setminus e'$ should be mapped onto one of the two ends of $f(e)$. The positions of the remaining $k-3$ vertices of $G$ can be chosen in $\binom{n-3}{k-3} (k-3)!$ different ways. In fact,
 in this way, some selections are counted more than once,
 namely those whose image contains both ends of $f(e)$.
We have $m$ choices for $e'$, therefore the number of
 non-$f$-exclusive copies of $G$ is smaller than
 $4m(k-3)\binom{n-3}{k-3} (k-3)!$.
Thus, as long as this number does not exceed $\binom{n}{k} k!$,
 every $f\in F_{n,0}$ generates an $f$-exclusive copy of $G$.
\end{proof}

The value of $g(G)$ for stars is known. 

\begin{theorem}[Alon, Caro \cite{AC}]\label{hkr1} For
     $r \ge  2$, we have $g(K_{1,r}) = 2r  $. 
\end{theorem}
\stepcounter{thm}

Concerning the function $g$, one can check that $g(K_2,1)=3$ and $g(K_2,0)=4$ hold; the difference comes from the fact that we need at least four vertices to be able to define a strong-shifting mapping. 
On the other hand, if $|f(e)\cap e|=1$ and $e_1,e_2,\dots,e_{t-1}$
 is an $f$-free copy of $(t-1)K_2$ with all $e_i$ being disjoint with $e$, then $e_1,\dots,e_{t-1},e$ is an $f$-free copy of $tK_2$. This shows that $g(tK_2,0) = g(tK_2,1)$ holds for any $t \ge 2$. So the next theorem determines all the values of $g(tK_2,0)$ as well.

\begin{thm}\label{tk2}
     $g(tK_2) =  2t$  for $t \ge 4$  and $g(tK_2) =2t +1$ for $t= 1,2,3$.  
\end{thm}

\begin{proof}
    The lower bound $2t$ is trivial as we need that many vertices to have a copy of $tK_2$. Also,  the number of $tK_2$ subgraphs in $K_{2t}$ is $a_{2t}=(2t-1)!! =  (2t-1) \cdot (2t-3)\cdot \ldots \cdot 3\cdot 1$. 

    %If a matching $M$ of size $t$ is not $f$-free, then there exists an edge $e$ such that $e,f(e)\in M$. Then we say that $e$ destroys $M$. Any edge $e\in E_{2t}$ destroys at most $a_{2t-4}$ matchings of size $t$. So the total number of destroyed matchings is at most $\binom{2t}{2}a_{2t-4}$. Rearranging $a_{2t}>\binom{2t}{2}a_{2t-4}$ yields the upper bound for $t\ge 4$.

    Let $f : E_{2t} \rightarrow E_{2t}$ be an edge mapping.
Say that an edge $e$ destroys a matching $M$ if $e, f(e)$ are both in $E(M)$. Any edge
$e \in E_{2t}$ destroys at most $a_{2t-4}$ matchings of size $t$. Thus, the number of destroyed matchings is at most 
$\binom{2t}{2}a_{2t-4}$ which is less than $a_{2t}$ when $t >3$. Hence, there
is a matching $M$ that is not destroyed by any edge which means that $M$ is $f$-free.
This proves $g(tK_2) \le 2t$ for $t \ge 4$.

For $t = 1$, clearly, we need three vertices to allow $f( e) \neq e$  and thus $g(K_2)= 3$.

For  $t = 2$    the only mapping $f\in F_{4,0}$ defined by $f(xy)=V(K_4)\setminus \{x,y\}$ shows $g(2K_2) \ge 5$.  To see the upper bound, observe that $K_5$   has precisely 15  pairs of $2K_2$.  Every edge $e$ can destroy at most one $2K_2$ so for any mapping $f\in F_{5,1}$ there remain at least five $f$-free $2K_2$.

Lastly, consider $t=3$. Observe that $K_6$ has 15 edges and 15 copies of $3K_2$. Introduce the bipartite graph $B$ with parts $E_6$ and the 15 copies of $3K_2$ in $K_6$ with $(e, M)\in E(B)$ if and only if $e\in M$. It is easy to see that $B$ is 3-regular and therefore has a perfect matching $\mathcal{M}$. So we can define a mapping $f\in F_{6,0}$ such that $f(e)=e^*$ with $e^*\neq e$ being one of the two edges in $M$ such that $(e, M)\in \mathcal{M}$. Then clearly, $f$ admits no free $3K_2$ and thus $g(3K_2)\ge 7$. The upper bound follows from the fact that an edge in $E_7$ can destroy at most three copies of $3K_2$ and that there are 105 copies of $3K_2$ in $K_7$ while $\binom{7}{2}\cdot 3<105$.
\end{proof}

%For another kind of upper bound on $g(G)$ and $g(G,0)$ we refer to Theorem~2.2 of \cite{C}.

\section{The appearance of $f$-exclusive copies of $G$: $w(G)$}

This section is devoted to the proof of Theorem~\ref{w} and also presents bounds on $w(K_k)$, as well as exact values of $w(K_{1,r})$ for infinitely many $r$.

\begin{proof}[\textbf{Proof of Theorem \ref{w}}]
First, we prove the upper bound on $w(G)$. Assume $k\ge 4$, and let
 $G$ be any graph with $k$ vertices and $m$ edges.
We fix a labeling $1,2,\dots,k$ on $V(G)$.
Also, we label the vertices of $K_n$ with $1,2,\dots,n$.
Then the number of labeled mappings from $V(G)$
 to $V(K_n)$ is equal to $n!/(n-k)!$.

Let $f\in F_{n,0}$ be any function.
Each $e\in E(K_n)$ can be extended in
 $2m\cdot (n-4)!/(n-k-2)!$ ways to a labeled
 subgraph $H\cong G$ such that $e\in E(H)$
 and $f(e)\cap V(H)=\emptyset$, while the
 total number of mappings from $G$ to $K_n$ that
 map an edge of $G$ to a given edge of $K_n$ is
 $2m\cdot (n-2)!/(n-k)!$.
This means that a fixed edge of $K_n$ is contained in
 $$
   2m \cdot
    \left[\, (n-2)(n-3)-(n-k)(n-k-1) \,\right] \cdot
   \frac{(n-4)!}{(n-k)!}
 $$
 non-$f$-exclusive images of $G$ under any chosen
  mapping $f\in F_{n,0}$.
If $\binom{n}{2}$ times this number is smaller than the
 number $n!/(n-k)!$ of labeled subgraphs $G$, then
 an $f$-exclusive copy of $G$ must occur for any $f$.
To avoid that this happens, we obtain the condition
 $$
   (n-2)(n-3) \leq m \cdot
    \left[\, (2k-4) n - (k^2 + k - 6) \,\right] .
 $$
Rearrangement yields
 $$
   n^2 - (2mk - 4m +5) n + (k+3)(k-2)m + 6 \le 0.
 $$
However, by plugging $n=2mk-4m+2$ into the left-hand side we get
$$
  (k+3)(k-2)m + 6 - 3(2mk - 4m +2) =
   m\cdot (k-2)(k-3)
$$
 which is positive for all $k\ge 4$, and therefore
 $w(G)\le 2mk-4m+2$ follows.

 Next, we prove the lower bound on $w(G)$. We take a mapping $f\in F_{n,0}$ at random, assuming that the
 choices $f(e)$ are totally independent for all $e$.
Let $G$ be a given graph with $k$ vertices and $m$ edges.
Then the number of vertex-labeled copies of $G$ in $K_n$ is
 $n!/(n-k)!<n^k$.
For each of the $m$ edges $e\in E(G)$ we have
 $\mathbb{P}(e\cap V(G)=\emptyset) =
  \binom{n-k}{2} / \binom{n-2}{2} < (1 - \frac{k-2}{n-2})^2$.
Since the choice of $f(e)$ is independent on the edges, the
 upper bound
  $$
    \mathbb{E}(\# \ \mbox{\rm of} \ f\mbox{\rm -exclusive} \
     G\subset K_n) < n^k \left(1 - \frac{k-2}{n-2}
     \right)^{2m}
    <\exp(k \ln n - 2m (k-2)/(n-2) )
  $$
 follows for the expected number of $f$-exclusive copies of $G$.
Consequently, as long as the exponent is negative, there exists
 an $f$ without $f$-exclusive $G$, which means $n<w(G)$.
\end{proof}

For complete graphs, we can obtain a tighter upper bound.
We also state a lower bound that is somewhat weaker but gives an
 explicit formula avoiding the $o(1)$ term.

\begin{thm}\label{wkk}
For any $k\ge 4$, we have $\frac{k^2}{3\log k}\le w(K_k)\le k(k-1)(k-2)+4- \lfloor\frac{k}{2}\rfloor$.
\end{thm}

\begin{proof}
To see the upper bound, consider an $f\in F_{n,0}$. Observe that for any edge $e\in E_n$ there are $2\binom{n-4}{k-3}+\binom{n-4}{k-4}$ $k$-cliques $K$ with $e\subset K$, $f(e)\cap K\neq \emptyset$. So if no $f$-exclusive $k$-cliques exist, then $\binom{n}{2}(2\binom{n-4}{k-3}+\binom{n-4}{k-4})\ge \binom{n}{k}$ holds. This is equivalent to 
\[
\frac{k(k-1)(k-2)}{(n-2)(n-3)}\left((n-k)+\frac{k-3}{2}\right)\ge 1.
\]
Plugging $n=k(k-1)(k-2)+4- \lfloor\frac{k}{2}\rfloor$, we obtain $(n-4+\lfloor\frac{k}{2}\rfloor)(n-\frac{k+3}{2})< (n-2)(n-3)$, a contradiction if $k\ge 3$.

To see the lower bound, we consider a mapping $f$ where $f(e)$ is picked uniformly at random from $\binom{[n]\setminus e}{2}$, for every edge $e\in E_n$. The expected number of $f$-exclusive $k$-cliques is $$\binom{n}{k}\left[\frac{\binom{n-k}{2}}{\binom{n}{2}}\right]^{\binom{k}{2}}\le \left(\frac{en}{k}\right)^k\left(\frac{n-k}{n}\right)^{2\binom{k}{2}}\le \left(\frac{en}{k}\right)^ke^{-\frac{k^2(k-1)}{n}}.$$
If $n\le \frac{k^2}{3\log k}$, then this is at most $e^{k+k\log k-k\log\log k-3(k-1)\log k}$, which is smaller than 1 if $k\ge 4$. Thus there exists an $f$ such that the number of $f$-exclusive $k$-cliques is 0.
\end{proof}
\noindent\textbf{Remark.} Let us also note that a recent result of Gishboliner, Jin and Sudakov~\cite{GJS}, explicitly motivated by our upper bounds in Theorem~\ref{w} and Theorem~\ref{wkk}, implies $w(G)=O(e(G))$ for graphs without isolated vertices, improving the general upper bound in broad ranges.

Our last result concerns $w(K_{1,r})$.

\begin{thm}\label{wstar}
\begin{enumerate}[(i)]
    \item For every $r\ge 2$, we have $w(K_{1,r})\le 5r-3$.
    \item The upper bound is sharp for every $r=(11^s+4)/5$ with $s\ge 1$. In particular, $w(K_{1,r})=5r-3$ for infinitely many values of $r$.
\end{enumerate}
\end{thm}

\begin{proof}
For the upper bound, let $f\in F_{5r-3,0}$. We apply Lemma~\ref{excluk1r}(i) with $G=H=K_{5r-3}$. Since $r\ge 2$, the graph $G=K_{5r-3}$ has at least four vertices, so no edge of $G$ is incident to all other edges; and we also have $\Delta(G)=5r-4$. Hence Lemma~\ref{excluk1r}(i) gives an $f$-exclusive copy of $K_{1,r}$ in $K_{5r-3}$. This proves $w(K_{1,r})\le 5r-3$.

For the lower bound for infinitely many values of $r$, let $q=11^s$ for some integer $s\ge 1$, and put $r=(q+4)/5$. Work in the finite field $\mathbb F_q$, and take $a=3\in \mathbb F_{11}\subseteq \mathbb F_q$. Then $a^2+a-1=0$, and $a$ has multiplicative order $5$.

Identify the vertex set of $K_q$ with $\mathbb F_q$. For an unordered edge $xy$, define
\[
    f(xy)=\{x+a(y-x),\ x+a^2(y-x)\}.
\]
This is well-defined: if we interchange $x$ and $y$, then the two image vertices become $y+a(x-y)$ and $y+a^2(x-y)$, which are the same two vertices since $1-a=a^2$ and $1-a^2=a$, both consequences of $a^2+a-1=0$. The map is strong-shifted: neither image vertex is equal to $x$, because $x\neq y$ and $a,a^2\neq 0$, and neither is equal to $y$ by the interchange property just observed.

Fix a possible center $x$ of a star. A leaf has the form $y=x+t$ with $t\in \mathbb F_q\setminus\{0\}$. The image of the edge $xy$ corresponds to the two directions $at$ and $a^2t$. Thus, on the $q-1$ possible leaves, the obstruction graph for an $f$-exclusive star is the disjoint union of complete graphs on the cosets of the order-$5$ subgroup $\langle a\rangle$ of $\mathbb F_q^*$. Indeed, the nontrivial ratios inside such a coset are $a,a^2,a^{-1},a^{-2}$. Therefore an $f$-exclusive star centered at $x$ can contain at most one leaf from each of the $(q-1)/5$ cosets. Hence every $f$-exclusive star has at most $(q-1)/5=r-1$ leaves, and this mapping contains no $f$-exclusive copy of $K_{1,r}$ in $K_q$.

Consequently $w(K_{1,r})>q=5r-4$. Together with the upper bound $w(K_{1,r})\le 5r-3$, this gives $w(K_{1,r})=5r-3$ for every $r=(11^s+4)/5$ with $s\ge 1$.
\end{proof}

\section{Closing remarks}

There are lots of open problems with respect to each of the parameters addressed in the paper. In this final section, we enumerate the ones that we find the most important or that seem to be natural next steps of future research. 

Theorem \ref{mtkr} gives reasonable bounds for $m(T,K_r)$, where $T$ is any tree on $k$ vertices.

\begin{problem}
    Determine $m(T,K_r)$ exactly or improve upon the lower and upper bounds given in Theorem \ref{mtkr}.
\end{problem}

Theorem \ref{th:m-tree} gives reasonable bounds on $m(T,K_{1,r})$, where $T$ is any tree on $k$ vertices.

\begin{problem}
    Determine $m(T,K_{1,r})$ exactly or improve upon the upper and lower bounds given in Theorem \ref{th:m-tree} items 2-3-4 in case $2(r-1)$ is not divisible by $k-1$.
\end{problem}

Theorem \ref{*tk} gives lower and upper bounds for $m^*(T,K_{1,r})$ that are still far apart.

\begin{problem}
    Determine $m^*(T,K_{1,r})$ exactly or improve upon the upper and lower bounds given in Theorem \ref{*tk}.
\end{problem}

Theorems~\ref{degreethm} and~\ref{gGkm} give bounds on $g(G)$ and $g(G,0)$, initially studied in \cite{AC,C}.

\begin{problem}
    Improve upon the bounds given in Theorems \ref{degreethm} and \ref{gGkm} and in particular find further graphs $G$, apart from stars and matchings, for which $g(G)$ and $g(G,0)$ can be determined exactly.
\end{problem}

Theorems~\ref{w} and~\ref{wkk} give lower and upper bounds for $w(G)$ and $w(K_k)$, and the upper bounds were improved in \cite{GJS}. Still, there is a gap to close.

\begin{problem}
    Improve upon the lower or upper bounds given in Theorem \ref{w} and Theorem \ref{wkk}. 
\end{problem} 

Theorem~\ref{wstar} determines $w(K_{1,r})$ for infinitely many values of $r$, in particular $w(K_{1,3})=12$. The authors also computed that $w(K_{1,2})=6$ and $8\le w(2K_2)\le 9$.

\begin{problem}
    Determine the exact values of $w(K_{1,r})$ for the remaining values of $r\ge 4$ and determine $w(tK_2)$ for $t\ge 2$.
\end{problem}

While the function $h(n, d, G)$ already appeared and was discussed in \cite{AC}, \cite{C} and \cite{CPTV} and few exact values are known \cite{CPTV}, the function $s(n, G)$ is introduced in \cite{CPTV} for the first time and is challenging when $G$ is bipartite (see Theorem 18 and Proposition 21 in \cite{CPTV}).

\begin{problem}
    Give lower and upper bounds for $s(n,G)$ for various graphs; of particular interest is finding graphs for which $s(n, G)$ can be determined exactly. 
\end{problem}


\begin{thebibliography}{99}
%\bibitem{AHS}
%Abbott, H.L., Hanson, D., Sauer, N., Intersection theorems for systems of sets, Journal of Combinatorial Theory, Series A, 12 (3) (1972), 381--389.

\bibitem{AC}
Alon, N., Caro, Y., Extremal problems concerning transformations of the set of edges of the complete graph. European J. Combinatorics 7(1986), 93--104.

\bibitem{ACT89}
Alon, N., Caro, Y., Tuza, Zs., Sub-Ramsey numbers for arithmetic progressions. Graphs and Combinatorics, 5(1) (1989), 307--314.

\bibitem{AM06} 
Axenovich, M., Martin, R., Sub-Ramsey numbers for arithmetic progressions. Graphs and Combinatorics, 22(3) (2006), 297--309.

\bibitem{C}
Caro, Y., Extremal problems concerning transformations of the edges of the complete hypergraphs. Journal of Graph Theory 11(1) (1987), 25--37.

\bibitem{CPTV}
Caro, Y., Patk\'os, B., Tuza, Zs., Vizer, M., Edge mappings of graphs: Tur\'an type parameters, European Journal of Combinatorics, 127 (2025) 104140

\bibitem{CS}
Caro, Y., Sch\"onheim, J., On edge-mappings with fixed edges avoiding free triangles. Ars Combinatoria 25-A (1988), 159--164.

\bibitem{chenetal}
Chen, G., Gould, R. J., Pfender, F., Wei, B., Extremal graphs for intersecting cliques. Journal of Combinatorial Theory, Series B, 89(2) (2003), 159--171.

\bibitem{Ch}
Chvátal, V., Tree-complete graph ramsey numbers. J. Graph Theory, 1 (1977) 93--93.


%\bibitem{CH}
%Chetwynd, A. G., Hilton, A. J. (1989). 1-factorizing regular graphs of high degree—an improved bound. Discrete Mathematics, 75(1-3), 103-112.
%\bibitem{CH}
%Chvátal, V., Hanson, D. Degrees and matchings. Journal of Combinatorial Theory, Series B, 20(2)  (1976), 128-138.


%\bibitem{CNR}
%Cutler, J., Nir, J. D., Radcliffe, A. J., Supersaturation for subgraph counts. Graphs and Combinatorics, 38(3) (2022), 1--21.

\bibitem{CFS16}
Conlon, D., Fox, J., Sudakov, B., Short proofs of some extremal results II. Journal of Combinatorial Theory, Series B, 121 (2016), 173--196.

\bibitem{CFS20}
Conlon, D., Fox, J., Sudakov, B., Short proofs of some extremal results III. Random Structures and Algorithms, 57(4) (2020), 958--982.

\bibitem{CFS-x}
Conlon, D., Fox, J., Sudakov, B., Independent arithmetic progressions. \url{https://arxiv.org/abs/1901.05084}

\bibitem{Ed}
Edwards, C. S., A lower bound for the largest number of triangles with a common edge (unpublished manuscript), 1977.

\bibitem{E64}
Erd\H os, P., Extremal problems in graph theory. In Theory of graphs and its applications, Proc. Sympos. Smolenice (1964), 29--36.

\bibitem{E69}
Erd\H os, P., On the number of complete subgraphs and circuits contained in graphs. \v Casopis P\v est. Mat. 94 (1969), 290--296.

\bibitem{erdosetal}
Erd\H os, P., F\"uredi, Z., Gould, R. J., Gunderson, D. S., Extremal graphs for intersecting triangles. Journal of Combinatorial Theory, Series B, 64(1) (1995), 89--100.

\bibitem{EH58}
Erd\H os, P., Hajnal, A., On the structure of set-mappings. Acta Mathematica Academia Scientia Hungarica, 9 (1958), 111--131.

\bibitem{ERS}
Erd\H os, P., R\'enyi, A., S\'os, V. T., On a problem of graph theory. Stud. Sci. Math. Hungar., 1 (1966), 215--235.

%\bibitem{ESa}
%Erd\H os, P., Sachs, H., Regul\"are graphen gegebener Taillenweite mit minimaler Knotenzahl. Wiss. Z. Martin-Luther-Univ. Halle-Wittenberg Math.-Natur. Reihe, 12 (1963) 251-257.

\bibitem{ES}
Erd\H os, P., Simonovits, M., A limit theorem in graph theory. Studia Sci. Math. Hungar., 1 (1966), 51--57.

\bibitem{ES83}
Erd\H os, P., Simonovits, M., Supersaturated graphs and hypergraphs. Combinatorica, 3(2) (1983), 181--192.

%\bibitem{FMS}
%Ferber, A., McKinley, G., Samotij, W., Supersaturated sparse graphs and hypergraphs. Internat. Math. Research Notices, Vol. 2020, Issue 2, 378--402.

\bibitem{FHMZ}
Fox, J., Himwich, Z., Mani, N., Zhou, Y., A note on directed analogues of the Sidorenko and forcing conjectures. \url{https://arxiv.org/abs/2210.16971}

\bibitem{FG}
F\"uredi, Z., Gunderson, D. S., Extremal numbers for odd cycles. Combinatorics, Probability \& Computing, 24(4) (2015), 641--645.

\bibitem{GJS}
Gishboliner, L., Jin, Z., Sudakov, B., Set mappings for general graphs. \url{arXiv:2601.00766}, 2026.

%\bibitem{GNV}
%Gerbner, D., Nagy, Z. L., Vizer, M., Unified approach to the generalized Tur\'an problem and supersaturation. Discrete Mathematics, 345(3) (2022), 112743.

\bibitem{HN}
Had\v ziivanov, N. G., Nikiforov, S. V., Solution of a problem of P. Erd\H os about the maximum number of triangles with a common edge in a graph (in Russian). C.R. Acad. Bulgare Sci., 32 (1979), 1315--1318.

\bibitem{kinai}
Hou, X., Qiu, Y., Liu, B., Extremal graph for Intersecting odd cycles. The Electronic Journal of Combinatorics, 23(2) (2016), P2.29

\bibitem{KZ}
Komatsu, T., Zhang, Y., Weighted Sylvester sums on the Frobenius set in more variables. \url{https://arxiv.org/abs/2101.04298}

\bibitem{LS}
Lov\'asz, L., Simonovits, M., On the number of complete subgraphs of a graph II. In: Studies in pure mathematics (1983), 459--495. Birkh\"auser, Basel.

\bibitem{MM}
Moon, J. W., Moser, L., On a problem of Tur\'an. Magyar Tudom\'anyos Akad\'emia Matematikai Kutat\'o Int\'ezeti K\"ozl\"ony, 7 (1962), 283--286.

\bibitem{R}
Reiher, C., The clique density theorem. Annals of Mathematics, 184(3) (2016), 683--707.

\bibitem{Sid}
Sidorenko, A. F., A correlation inequality for bipartite graphs. Graphs and Combinatorics 9 (1993), 201--204.

\bibitem{S}
Simonovits, M., Extremal graph problems with symmetrical extremal graphs. Additional chromatic conditions. Discrete Mathematics, 1 (1974), 349--376.

\bibitem{S84}
Simonovits, M., Extremal graph problems, degenerate extremal problems, and supersaturated graphs. In: Progress in Graph Theory (Waterloo, Ont., 1982), (1984), 419--437.

\bibitem{TT}
Tait, M., Timmons, C., Regular Tur\'an numbers of complete bipartite graphs. Discrete Mathematics, 344(10) (2021), 112531.

\end{thebibliography}
\end{document}